\documentclass[10pt,twoside]{article}
\usepackage{graphicx}
\usepackage{amsmath}
\usepackage{Latex-document}

\def\RR{\rm \hbox{I\kern-.2em\hbox{R}}}
\def\NN{\rm \hbox{I\kern-.2em\hbox{N}}}
\def\ZZ{\rm {{\rm Z}\kern-.28em{\rm Z}}}
\def\CC{\rm \hbox{C\kern -.5em {\raise .32ex \hbox{$\scriptscriptstyle |$}}\kern
-.22em{\raise .6ex \hbox{$\scriptscriptstyle |$}}\kern .4em}}

\def\<{\langle}
\def\>{\rangle}
\def\t{\tilde}

\def\e{\varepsilon}

\def\nl{\newline}

\def\cT{{\cal T}}
\def\cF{{\cal F}}
\def\Chi{\raise .3ex
\hbox{\large $\chi$}} 

\def\[{\Bigl [}
\def\]{\Bigr ]}
\def\({\Bigl (}
\def\){\Bigr )}
\def\[{\Bigl [}
\def\]{\Bigr ]}
\def\({\Bigl (}
\def\){\Bigr )}

\newcommand{\be}{\begin{equation}}
\newcommand{\ee}{\end{equation}}
\newcommand{\bea}{\begin{array}{lll}}
\newcommand{\eea}{\end{array}}
\newcommand{\bi}{\begin{itemize}}
\newcommand{\ei}{\end{itemize}}
\newcommand{\iref}[1]{(\ref{#1})}
\def\volumeno{I}

\title{\bf  Adaptive Methods for PDE's \vskip -2mm
Wavelets or Mesh Refinement? \vskip 6mm}
\author{Albert Cohen\vspace*{-0.5cm}\thanks{Laboratoire Jacques-Louis
Lions, Universit\'e Pierre et Marie Curie, Paris, France. E-mail:
cohen@ann.jussieu.fr }}
\date{\vspace{-8mm}}

\begin{document}

\maketitle

\thispagestyle{first} \setcounter{page}{607}

\begin{abstract}

\vskip 3mm

Adaptive mesh refinement techniques are nowadays an established
and powerful tool for the numerical discretization  of PDE's. In
recent years, wavelet bases have been proposed as an alternative
to these techniques. The main motivation for the use of such bases
in this context is their good performances in data compression and
the approximation theoretic foundations which allow to analyze and
optimize these performances. We shall discuss these theoretical
foundations, as well as one of the approaches which has been
followed  in developing efficient adaptive wavelet solvers. We
shall also discuss the similarities and differences between
wavelet methods and adaptive mesh refinement.

\vskip 4.5mm


\end{abstract}

\vskip 12mm

\section*{1. Introduction}

\vskip-5mm \hspace{5mm}

Among those relevant phenomenons which are modelled by partial
differential or integral equations, countless are the instances
where the mathematical solutions exhibit {\it singularities}.
Perhaps the most classical examples are elliptic equations on
domains with re-entrant corners, or nonlinear hyperbolic systems
of conservation laws. While such singularities are sources of
obvious theoretical difficulties---classical solutions should be
abandonned to the profit of weak solutions---they are also an
obstacle to the convergence of numerical approximation methods, in
the sense that they deteriorate the rate of decay of the error
with respect to the size of the discrete problem : achieving a
prescribed accuracy will typically require finer resolution and
therefore heavier computational cost and memory storage, in
comparison to the approximation of smooth solutions. Let us remark
that singularities often have a physical relevance : they
represent the concentration of stress in elasticity, boundary
layers in viscous fluid flows, shock waves in gas dynamics... It
is therefore a legitimous requirement that they should be
accurately resolved by the numerical method.

In this context, the use of adaptive methods,
appears as a natural solution to improve the approximation at a reasonable
computational cost. Here, the word {\it adaptivity} has a twofold meaning :
(i) the discretization is allowed to be refined only locally, in particular near the singularities of the
solution, and (ii) the resolution algorithm uses information
gained during a given stage of the
computation in order to derive a new refined discretization for the next stage.
The most typical example is {\it adaptive mesh refinement}
based on {\it a-posteriori} error estimates in the finite element context.
While these methods have proved to be computationally successful,
the theory describing their advantages over their non-adaptive counterpart
is far from being complete. In particular, the {\it rate of convergence}
of the adaptive algorithm, which describes the trade-off
between the accuracy and
complexity of the approximation, is not clearly understood.

In recent years,
wavelet bases have been proposed as an alternative to
adaptive mesh refinement, motivated by their good performances
in data (more specifically image) compression.
In wavelet-based adaptive schemes,
the set of basis functions which describe the approximate solution is
updated at each stage of the computation.
Intuitively, the selection of the appropriate basis functions plays a similar role
as the selection of the mesh points in adaptive finite element methods,
and one could therefore expect similar performances from both approaches.
On a more rigorous level, a specific feature of the wavelet approach
is the emergence of a sound theoretical setting which
allows to tackle fundational questions such as the rate of convergence
of the adaptive method.

The goal of this paper is to give some elements of
comparison between adaptive wavelet and mesh refinement methods
from this perspective. We shall
first describe in \S 2 a general setting which leads us in \S 3 to
a first comparison between wavelets and adaptive finite elements from the point
of view of approximation theory. We discuss in \S 4 the relation between
these results and adaptive algorithms for PDE's. After recalling  in \S 5
the classical approach in the finite element context, we present in \S 6
an adaptive wavelet strategy which has been applied to various problems,
and discuss its fundational specificities. Finally, we shall
conclude in \S 7 by pointing out some intrinsic shortcoming
of wavelet-based adaptive methods.

\section*{2. A general framework}

\vskip-5mm \hspace{5mm}

Approximation theory is the branch of mathematics which studies the
process of approximating general functions by simple functions
such as polynomials, finite elements or Fourier series. It plays
therefore a central role in the accuracy analysis of numerical
methods. Numerous problems of approximation theory have in common
the following general setting : we are given a
family of  subspaces $(S_N)_{N\geq 0}$
of a normed space $X$, and for $f\in X$, we consider the {\it best approximation error}
\be
\label{infimum}
\sigma_N(f):=\inf_{g\in S_N}\|f-g\|_{X}.
\ee
Typically, $N$ represents the number of parameters which are needed
to describe an element in $S_N$, and in most cases of interest,
$\sigma_N(f)$ goes to zero as this number tends to infinity.
If in addition $\sigma_N(f) \leq CN^{-s}$ for some $s>0$, we say that $f$ is
approximated at rate $s$.

Given such a setting, the central problem of approximation theory
is to {\it characterize} by some analytic (typically smoothness)
condition those functions $f$ which are approximated at some
prescribed rate $s>0$. Another important problem is how to design
simple approximation procedures $f\mapsto f_N\in \Sigma_N$ which
avoid solving the minimization problem \iref{infimum}, while
remaining near optimal in the sense that \be \|f-f_N\|_X \leq
C\sigma_N(f), \ee for some constant $C$ independent of $N$ and
$f$.

As an example, consider approximation
by finite element spaces $V_h$ defined from regular conforming partitions
$\cT_h$ of a domain $\Omega\subset \RR^d$ into simplices
with uniform mesh size $h$. The approximation
theory for such spaces is quite classical, see e.g. \cite{Cia},
and can be summarized in the following way.
If $W^{t,p}$ denotes the classical Sobolev space, consisting of
those functions $f\in L^p$ such that $D^\alpha f\in L^p$ for $|\alpha|\leq t$, we typically
have
\be
f\in W^{t+r,p}\Rightarrow \inf_{g\in V_h}\|f-g\|_{W^{t,p}}\leq Ch^{r}
\ee
provided that $V_h$ is contained in $W^{t,p}$ and that
$V_h$ has approximation order larger than $t+r$, i.e.
contains all polynomials of degree strictly less than
$t+r$. Such classical results also hold for fractional smoothness.
We can express them in terms of the number of parameters,
remarking that $N:={\rm dim}(V_h)\sim h^{-d}$, so that if we set
$X=W^{t,p}$ and $S_N:=V_h$ with $h:=N^{-1/d}$, we have obtained
\be
f\in W^{t+r,p}\Rightarrow  \sigma_N(f)\leq C N^{-r/d}.
\label{linearerror}
\ee
We have thus identified an analytic condition which ensures the
rate $s=r/d$. Note that this is not a characterization (we only have an implication),
yet a deeper analysis shows that an ``if and only if'' result holds if we slightly
modify the notion of Sobolev smoothness (using Besov classes, see \cite{Co}).
In summary, the rate of approximation in $W^{s,p}$ is
governed by the approximation order of the $V_h$ spaces, the dimension $d$
and the level of smoothness of $f$ measured in $L^p$. Let us finally
remark that near-optimal approximation procedures can be obtained
if we can find a sequence of finite element projectors $P_N: X \mapsto S_N$ such that
$\|P_N\|_{X\to X}\leq K$ with $K$ independent of $N$ : in this case, we simply take
$f_N=P_Nf$ and remark that $\|f-f_N\|_X \leq (1+K)\sigma_N(f)$.

In the following we shall address the same
questions in the cases of adaptive finite element and wavelet
approximation. As we shall see, a specific feature to such cases is that
the spaces $\Sigma_N$ are not linear vector spaces.

\section*{3. Adaptive finite elements and wavelets}

\vskip-5mm \hspace{5mm}

In the adaptive finite element setting, the number of parameters
$N$ is proportional to the number of triangles, but for a given
budget $N$ the partition $\cT$ and the finite element space
$V_\cT$ are allowed to be locally refined in a way which depends
on the function $f$ to be approximated. It is therefore natural to
define the approximation spaces $S_N$ as \be S_N:=\cup_{
\#(\cT)\leq N}  V_\cT. \ee It should be well understood that the
$S_N$ are not linear vector spaces (the sum of two elements does
not in general fall in $S_N$ when their triangulation do not
match) but any $g\in S_N$ is still described by ${\cal O}(N)$
parameters, which encode both its triangulation $\cT$ and its
coordinates in $V_\cT$. The requirement of adaptivity has thus led
us to the concept of  {\it nonlinear approximation}.

Wavelet bases offer another track toward nonlinear adaptive
approximation. The simplest prototype of a wavelet basis is the
{\it Haar system}. Let us describe this system in the case of
expanding a function $f$ defined on $[0,1]$ : the first component
in this expansion is simply the average of $f$, i.e. the
orthogonal projection $\<f,e_0\>e_0$ onto the function
$e_0=\Chi_{[0,1]}$. The approximation is then refined into the
average of $f$ on the two half intervals of equal size. This
refinement amounts in adding the orthogonal projection
$\<f,e_1\>e_1$ onto the function
$e_1=\Chi_{[0,1/2]}-\Chi_{[1/2,1]}$. Iterating this refinement
process, we see that the next components have the same form as
$e_1$ up to a change of scale : at refinement level $j$, we are
adding the projection onto the functions \be
\psi_{j,k}(x)=2^{j/2}\psi(2^jx-k),\;\; k=0,\cdots,2^{j}-1, \ee
where $\psi=e_1$. Since all these functions are orthogonal to the
previous ones, letting $j$ go to $+\infty$, we obtain the
expansion of $f$ into an orthonormal basis of $L^2([0,1])$ \be
f=\sum_{\lambda\in\nabla}f_\lambda\psi_\lambda, \ee with
$f_\lambda:=\<f,\psi_\lambda\>$. In the above notation $\lambda$
concatenates the scale and space parameters $j$ and $k$, and
$\nabla$ is the set of all indices (including also the first
function $e_0$). In order to keep track of the scale $j$
corresponding to an index $\lambda=(j,k)$ we shall use the
notation $|\lambda|=j$. More general wavelet systems in one or
several space dimension are built from similar nested
approximation processes, involving e.g. spline functions or finite
elements in place of piecewise constant functions (see \cite{Dau}
or \cite{Co} for a general presentation).

This brief description suggests that a natural
construction of adaptive wavelet approximations is obtained
by using only a limited set of indices $\lambda$ as the scale $|\lambda|$ grows,
which depends on the function to be approximated and typically
corresponds to those wavelets whose supports are close to its singularities. It is therefore natural to define
the approximation
spaces $S_N$ as the set of all $N$ terms combinations
\be
S_N:=\{\sum_{\lambda\in \Lambda}d_\lambda\psi_\lambda\;\; ;\;\; \#(\Lambda)\leq N\}.
\ee
Again this is obviously not a linear space, since we allow to approximate a function by choosing
the best $N$ terms which differ from one function to another. Note that
we still do have $S_N+S_N=S_{2N}$.

Both adaptive finite element and wavelet framework have obvious
similiarities. However, the answer to the two basic questions
raised in the previous section---what are the properties of $f$
which govern the decay of $\sigma_N(f)$ and how to compute in a
simple way a near optimal approximation of $f$ in $\Sigma_N$ ---is
only fully understood in the wavelet framework. Concerning the
first question, a striking result  by DeVore and his collaborators
\cite{De} is the following : with $X:=W^{t,p}$, best $N$-term
wavelet approximation satisfies \be \label{nonlinearerror} f\in
W^{t+r,q}\Rightarrow  \sigma_N(f)\leq C N^{-r/d}, \ee with $q$ and
$r$ connected by the relation $1/q=1/p+r/d$, assuming that the
multiresolution approximation spaces associated to the wavelet
basis are in $W^{t,p}$ and contain the polynomials of degree
strictly less than $t+r$.

Such an estimate should be compared with the linear estimate
\iref{linearerror} : the same convergence rate is governed by a
much weaker smoothness assumption on $f$ since $q<p$ (as in the
linear case, an ``iff and only if'' result can be obtained up to
slight technical modifications in the statement of
\iref{nonlinearerror}). This result gives a precise mathematical
meaning to the spatial adaptation properties of best $N$-term
wavelet approximation: a function $f$ having isolated
discontinuity, has usually a smaller amount of smoothness $s+t$
when measured in $L^p$ than when measured in $L^q$ with
$1/q=1/p+t/d$, and therefore $\sigma_N(f)$ might decrease
significantly faster than $\e_N(f)$.

The answer to the second question is given by a result due to
Temlyakov : if $f=\sum_{\lambda\in\nabla}d_\lambda\psi_\lambda$, and if we measure the
approximation error in $X=W^{t,p}$, a
near optimal strategy when $1<p<\infty$ consists in
the thresholding procedure which retains the $N$ largest contributions
$\|d_\lambda\psi_\lambda\|_{X}$ : if
$\Lambda_N$ is the corresponding set of indices, one can prove
that there exists $C>0$ independent of
$N$ and $f$ such that
\be
\|f-\sum_{\lambda\in\Lambda_N}d_\lambda\psi_\lambda\|_{X}\leq C \sigma_N(f).
\ee
This fact is obvious when $X=L^2$ using the orthonormal basis property.
It is a remarkable property of wavelet bases that it also holds for
more general function spaces. In summary, thresholding plays for best $N$-term wavelet
approximation an analogous role
as projection for linear finite element approximation.

In the adaptive finite element framework, a similar theory is
far from being complete. Partial answers to the basic questions
are available if one chooses to
consider adaptive partitions with {\it shape constraints} in terms of
a uniform bound on the aspect ratio of the elements
\be
\max_{K\in\cT} \([{\rm Diam}(K)]^d/ {\rm vol}(K) \)\leq C.
\ee
Such a restriction means that the local refinement is isotropic,
in a similar way to wavelets. In such a case, we therefore
expect a rate of approximation similar to \iref{nonlinearerror}.
Such a result is not available, yet the following can be proved \cite{Co} for Lagrange
finite elements of degree $m$ :
if for any given tolerance $\e>0$, one is able to build a partition
$\cT=\cT(\e)$ of cardinality $N=N(\e)$ such that on each $K\in \cT$ the local error of approximation
by polynomials satisfies
\be
\frac \e 2 \leq \inf_{p\in\Pi_m} \|f-p\|_{W^{t,p}(K)} \leq \e,
\ee
then we can build global approximants $f_N\in V_\cT \subset S_N$ such that
\be
f\in W^{t+r,q} \Rightarrow \|f-f_N\|_{X}\leq CN^{-r/d},
\ee
with $q$ and $r$ connected by the relation $1/q=1/p+r/d$ and assuming $s+t<m$.
The effective construction of $\cT(\e)$ is not always feasible, in particular
due to the conformity constraints on the partition which does not allow
to connect very coarse and very fine elements without intermediate grading.
However, this result shows that from an intuitive point of view, the adaptive finite
element counterpart to wavelet thresholding amounts in
{\it equilibrating the local error} over the partition.
One can actually use these ideas in order to obtain the estimate \iref{nonlinearerror}
for adaptive finite elements under the more restrictive assumption that $1/q<1/p+r/d$.
Let us finally mention that the approximation theory for
adaptive finite elements without shape constraints is an open problem.

\section*{4. Nonlinear approximation and PDE's}

\vskip-5mm \hspace{5mm}

Nonlinear approximation theory has opened new lines
of research on the theory of PDE's and their numerical discretization.
On the one hand, it is worth revisiting the regularity theory of
certain PDE's for which the solutions develop singularities
but might possess significantly higher smoothness in the
scale of function spaces which govern the rate of nonlinear approximation
in a given norm than in the scale which
govern the rate of linear approximation in the same norm.
Results of this type have been proved in particular
for  elliptic problems on nonsmooth domains \cite{DD}
and for scalar 1D conservation laws \cite{DeLu}.
These results show that if $u$ is the solution of such equations,
the rate of decay of $\sigma_N(u)$ is significantly
higher for best $N$-term approximation
than for the projection on uniform finite element spaces,
therefore advocating for the use of
adaptive discretizations of such PDE's.

On the other hand these results also provide with an ideal {\it
benchmark} for adaptive discretizations of the equation, since
$\sigma_N(u)$ represents the best accuracy which can be achieved
by $N$ parameters. In the wavelet case these parameters are
typically the $N$ largest coefficients of the exact solution $u$.
However, in the practice of solving a PDE, these coefficients are
not known, and neither is the set $\Lambda$ corresponding to the
indices of the $N$ largest contributions
$\|d_\lambda\psi_\lambda\|$. It is therefore needed to develop
appropriate {\it adaptive resolution strategies} as a substitute
to the thresholding procedure. Such strategies aim at detecting
the indices of the largest coefficients of the solutions and to
compute them accurately, in a similar way that adaptive mesh
refinement strategies aim at contructing the optimal mesh for
finite element approximation. In both contexts, we could hope for
an algorithm which builds approximations $u_N\in \Sigma_N$ such
that $\|u-u_N\|_X$ is bounded up to a fixed multiplicative
constant by $\sigma_N(u)$ for a given norm of interest, but this
requirement is so far out of reach. A more reasonable goal is that
the adaptive strategy exhibits the optimal rate of approximation :
if $\sigma_N(u) \leq C N^{-s}$ for some $s>0$, then $\|u-u_N\|_X
\leq CN^{-s}$ up to a change in the constant. Another requirement
is that the adaptive algorithm should be {\it scalable}, i.e. the
number of elementary operations in order to compute $u_N$ remains
proportional to $N$. Let us finally remark that the norm $\|\cdot
\|_X$ for which error estimates can be obtained is often dictated
by the nature of the equation (for example $X=H^1$ in the case of
a second order elliptic problem) and that additional difficulties
can be expected if one searches for estimates in a different norm.

\section*{5. The classical approach}

\vskip-5mm \hspace{5mm}

The {\em classical} approach to numerically solving linear and nonlinear
partial differential or integral equations $\cF(u)=0$ by the finite element method
is typically concerned with the following issues :
\begin{itemize}
\item[(c1)]
Well-posedness of the equation, i.e. existence, uniqueness and stability of the solution.
\item[(c2)]
Discretization into a finite element problem $\cF_\cT(u_\cT)=0$ by the Galerkin method with
$u_\cT\in V_\cT$, analysis of well-posedness and of the approximation error $\|u-u_\cT\|_X$.
\item[(c3)]
Numerical resolution of the finite dimensional system.
\item[(c4)]
Mesh refinement based on a-posteriori error estimators in the case of adaptive
finite element methods.
\end{itemize}
Several difficulties are associated to
each of these steps. First of all, note that the well-posedness of the finite element
problem is in general {\em not} a consequence
of the well-posedness of the continuous problem. Typical examples
even in the linear case are {\em saddle point
problems}. For such problems, it is well known that, for Galerkin discretizations to be
stable, the finite element spaces for the different solution components
have to satisfy certain compatibility conditions (LBB or Inf-Sup condition),
which are also crucial in the derivation of optimal error estimates.
Thus the discrete problem does not necessarily inherit
the ``nice properties'' of the original infinite dimensional
problem. Concerning the numerical resolution of the discrete system, a typical
source of trouble is its possible {\em ill-conditioning}, which interferes with the typical need
to resort on iterative solvers in high dimension.
An additional difficulty occuring in the case of
integral equations is the manipulation of
matrices which are densely populated.

Finally, let us elaborate more on the adaptivity step. Since more
than two decades, the understanding and practical realization of
adaptive refinement schemes in a finite element context has been
documented in numerous publications \cite{BR,BM,BW,EEHJ,Ve}. Key
ingredients in most adaptive algorithms are {\em a-posteriori
error estimators} which are typically derived from the current
residual $\cF(u_\cT)$ : in the case where the Frechet derivative
$D\cF(u)$ is an isomorphism between Banach function spaces $X$ to
$Y$, one can hope to estimate the error $\|u-u_\cT\|_X$ by the
evaluation of $\|\cF(u_\cT)\|_Y$. The rule of thumb is then to
decompose $\|\cF(u_\cT)\|_Y$ into computable local error
indicators $\eta_K$ which aim to describe as accurately as
possible the local error on each element $K\in \cT$. In the case
of elliptic problems, these indicators typically consist of local
residuals and other quantities such as jumps of derivatives across
the interface between adjacent elements. A typical refinement
algorithm will subdivide those elements $K$ for which the error
indicator $\eta_K$ is larger than a prescribed tolerance $\e$
resulting in a new mesh $\t \cT$. Note that this strategy is
theoretically in accordance with our remarks in \S 3 on adaptive
finite element approximation, since it tends to {\it equilibrate}
the local error. Two other frequently used strategies consist in
refining a fixed proportion of the elements corresponding to the
largest $\eta_K$, or the smallest number of elements $K$ for which
the $\eta_K$ contribute  to the global error up to a fixed
proportion. It is therefore hoped that the iteration of this
process from an initial mesh $\cT_0$ will produce optimal meshes
$(\cT_n)_{n\geq 0}$ in the sense that the associated solutions
$u_n:=u_{\cT_n}\in V_{\cT_n}$ converge to $u$ at the optimal rate
: \be \label{idealfem} \sigma_N(u)\leq CN^{-r} \Rightarrow
\|u-u_n\|_X \leq C [\#(\cT_n)]^{-r}, \ee up to a change in the
constant $C$. Unfortunately, severe obstructions appear when
trying to prove \iref{idealfem} even in the simplest model
situations. One of them is that $\eta_K$ is in general not an
estimate by above of the local error, reducing the chances to
derive the optimal rate. For most adaptive refinement algorithms,
the theoretical situation is actually even worse in the sense that
it cannot even be proved that the refinement step actually results
in a reduction of the error by a fixed amount and that $u_n$
converges to $u$ as $n$ grows. Only recently \cite{Do,MNS} have
proof of convergence appeared for certain type of adaptive finite
element methods, yet without convergence rate and therefore no
guaranteed advantage over their non-adaptive counterparts.

\section*{6. A new paradigm}

\vskip-5mm \hspace{5mm}

Wavelet methods vary from finite element method in that they can
be viewed as solving systems that are finite sections of one fixed
infinite dimensional system corresponding to the discretization of
the equation in the full basis. This observation has led to a {\em
new paradigm} which has been explored in \cite{CDD2} for {\em
linear variational problems}. It aims at closely intertwining the
analysis---discretization---solution process. The basic steps
there read as follows :
\begin{itemize}
\item[(n1)]
Well-posedness of the variational problem.
\item[(n2)]
Discretization into an
{\em equivalent} infinite dimensional problem which is
well posed in $\ell^2$.
\item[(n3)]
Devise an iterative scheme for the $\ell_2$-problem
that exhibits a fixed error reduction per iteration step.
\item[(n4)]
Numerical realization of the iterative scheme by means of
an {\em adaptive application} of the involved infinite dimensional
operators within some dynamically updated accuracy tolerances.
\end{itemize}
Thus the starting point (n1) is the same.
The main difference is that one aims at staying as long as
possible with the infinite dimensional problem. Only at the very end, when it comes
to applying the operators in the ideal iteration scheme (n4),
one enters the finite dimensional realm. However, the finite
number of degrees of freedom is determined at each stage by the
adaptive application of the operator, so that at no stage
any specific trial space is fixed.

The simplest example is provided by the Poisson equation $-\Delta u=f$
on a domain $\Omega$ with homogeneous boundary conditions, for which the
variational formulation in $X=H^1_0$ reads : find $u\in X$ such that
\be
a(u,v)=L(v),\;\; \mbox{for all} \;\; v\in X,
\ee
with $a(u,v):=\int_\Omega \nabla u\nabla v$ and $L(v):=\int_\Omega f v$. The well-posedness
for a data $f\in X'=H^{-1}$ is ensured by the Lax-Milgram lemma. In the analysis
of the wavelet discretization of this problem, we shall invoke the fact that
wavelet bases provide norm equivalence for Sobolev spaces in terms
of weighted $\ell^2$ norms of the coefficients : if $u=\sum_{\lambda}u_\lambda\psi_\lambda$, one
has
\be
\|u\|_{H^s}^2 \sim \sum_\lambda \|u_\lambda\psi_\lambda\|_{H^s}^2
\sim \sum_{\lambda} 2^{2s|\lambda|}|u_\lambda|^2.
\ee
We refer to \cite{Co} and \cite{Da} for the general mechanism allowing to derive these
equivalences, in particular for Sobolev spaces on domains with boundary conditions such as $H^1_0$.
Therefore, if we renormalize our system in such a way that $\|\psi_\lambda\|_X=1$, we obtain
the norm equivalence
\be
\|u\|_X^2 \sim \|U\|^2,
\ee
where $U:=(u_\lambda)_{\lambda\in\nabla}$ and $\|\cdot\|$ denotes the $\ell^2$ norm.
By duality, one also easily obtains
\be
\|f\|_{X'}^2 \sim \|F\|^2,
\ee
with $F:=(\<f,\psi_\lambda\>)_{\lambda\in\nabla}$. The equivalent $\ell^2$ system is thus given by
\be
AU=F,
\ee
where $A(\lambda,\mu)=a(\psi_\lambda,\psi_\mu)$ is a symmetric positive definite matrix
which is continuous and coercive in $\ell^2$. In this case, a converging infinite dimensional
algorithm can simply be obtained by the Richardson iteration
\be
U^n:=U^{n-1}+\tau (F-AU^{n-1})
\ee
with $0<\tau <2[\lambda_{\max}(A)]^{-1}$ and $U^0=0$, which guarantees
the reduction rate $\|U-U^n\| \leq \rho \|U-U^{n-1}\|$ with
$\rho= \max \{1-\tau\lambda_{\min}(A),\tau\lambda_{\max}(A)-1\}$.
Note that renormalizing the wavelet system plays the role of a multiscale preconditioning, similar to multigrid
yet operated at the infinite dimensional level.

At this stage, one enters finite dimensional adaptive computation
by modifying the Richardson iteration up to a prescribed tolerance
according to \be \label{modifiedit} U^n :=U^{n-1}+\tau ({\bf
COARSE}(F,\e)-{\bf APPROX}(AU^{n-1},\e)) \ee where $\|F-{\bf
COARSE}(F,\e)\|\leq \e$ and $\|AU-{\bf APPROX}(AU^{n-1},\e)\|\leq
\e$, and the $U^n$ are now finite dimensional vector supported by
adaptive sets of indices $\Lambda_n$. The procedure {\bf COARSE},
which simply corresponds to thresholding the data vector $F$ at a
level corresponding to accuracy $\e$, can be practically achieved
without the full knowledge of $F$ by using some a-priori bounds on
the size of the coefficients $\<f,\psi_\lambda\>$, exploiting the
local smoothness of $f$ and the oscillation properties of the
wavelets. The procedure {\bf APPROX} deserves more attention. In
order to limitate the spreading effect of the matrix $A$, one
invokes its {\it compressibility} properties, namely the possibily
to truncate it into a matrix $A_N$ with $N$ non-zero entries per
rows and columns in such a way that \be \|A-A_N\|_{\ell^2\to
\ell^2}Ê\leq CN^{-s}. \ee The rate of compressibility $s$ depends
on the available a-priori estimates on the off-diagonal entries
$A(\lambda,\mu):=\int_{\Omega}\nabla \psi_\lambda\nabla\psi_\mu$
which are consequences of the smoothness and vanishing moment
properties of the wavelet system, see \cite{CDD1}. Once these
properties are established, a first possibility is thus to choose
\be \label{approx1} {\bf APPROX}(AU^{n-1},\e)=A_N U^{n-1} \ee with
$N$ large enough so that accuracy $\e$ is ensured. Clearly the
modified iteration \iref{modifiedit} satisfies $\|U-U^n\| \leq
\rho \|U-U^{n-1}\| + 2\tau\e$, and therefore ensures a fixed
reduction rate until the error is of the order $\frac{2\tau}{1 -
\rho} \e$, or until the residual $F-AU^n$ is of order $\frac{2
\tau \| A \|}{1 - \rho} \e$. A natural idea is therefore to update
dynamically the tolerance $\e$, which is first set to $1$ and
divided by $2$ each time the approximate residual ${\bf
COARSE}(F,\e)-{\bf APPROX}(AU^{n-1},\e)$ is below $[\frac{2 \tau
\| A \|}{1 - \rho} + 3] \e$ (which is ensured to happen after a
fixed number of steps).

We therefore obtain a converging adaptive strategy, so far
without information about the convergence rate. It turns out
that the optimal convergence rate can also be proved, with
a more careful tuning of the adaptive algorithm.
Two additional ingredients are involved in this tuning.

Firstly, the adaptive matrix vector multiplication {\bf APPROX} has to be designed
in a more elaborate way than \iref{approx1} which could have the effect of
inflating too much the sets $\Lambda_n$. Instead, one defines for a finite length vector $V$
\be
\label{approx2}
{\bf APPROX}(AV,\e)=\sum_{l=0}^jA_{2^{j-l}} [V_{2^l}-V_{2^{l-1}}]
\ee
where $V_N$ denotes the restriction of $V$ to its $N$ largest components
(with the notation $V_{1/2}=0$), and $j$ is the smallest positive integer such
that the residual $\sum_{l=0}^j\|A-A_{2^{j-l}}\| \| V_{2^l}-V_{2^{l-1}}\| + \|A\|  \|V-V_{2^j}\|$
is less than $\e$. In this procedure, the spreading of the operator is
more important on the largest coefficients which are less in number,
resulting in a significant gain in the complexity of the outcome.

Secondly, additional coarsening steps are needed in order to
further limitate the spreading of the sets $\Lambda_n$ and preserve the
optimal rate of convergence. More precisely,
the procedure {\bf COARSE} is applied to $U^n$
with a tolerance proportional to $\e$, for those
$n$ such that $\e$ will be updated at the next iteration.

With such additional ingredients, it was proved in \cite{CDD2}
that the error has the optimal rate of decay in the sense that \be
\label{idealwave} \sigma_N(u)\leq CN^{-s} \Rightarrow \|u-u_n\|_X
\sim \|U-U^n\| \leq C [\#(\Lambda_n)]^{-s}, \ee and that moreover,
the computational cost of producing $u_n$ remains proportional to
$\#(\Lambda_n)$. It is interesting to note that this strategy
extends to non-elliptic problems such as saddle-points problems,
without the need for compatibility conditions, since one inherits
the well-posedness of the continuous problem which allows to
obtain a converging infinite dimensional iteration, such the Uzawa
algorithm or a gradient descent applied to the least-square system
(see also \cite{CDD2,DDU}). The extension to nonlinear variational
problems, based on infinite dimensional relaxation or Newton
iterations, has also been considered in \cite{CDD3}. It requires a
specific procedure for the application of the nonlinear operator
in the wavelet coefficients domain which generalizes
\iref{approx1}. It should also be mentioned that matrix
compressibility also applies in the case of integral operators
which have quasi-sparse wavelet discretizations. Therefore several
of the obstructions from the classical approach---conditioning,
compatibility, dense matrices---have disappeared in the wavelet
approach.

Let us finally mention that the coarsening steps are not really needed
in the practical implementations of the adaptive wavelet method
(for those problems which have been considered so far) which still does exhibit
optimal convergence rate. However, we do not know how to prove
\iref{idealwave} without these coarsening steps. There seems to be a
similar situation in the finite element context :
it has recently been proved in \cite{BDD} that \iref{idealfem}
can be achieved by an adaptive mesh refinement algorithm which incorporates
coarsening steps, while these steps are not needed in practice.

\section*{7.  Conclusions and shortcomings}

\vskip-5mm \hspace{5mm}

There exist other approaches for the development of efficient
wavelet-based adaptive schemes. In particular, an substantial
research activity has recently been devoted to {\it
multiresolution adaptive processing} techniques, following the
line of idea introduced in \cite{H1,Ha}. In this approach, one
starts from a classical and reliable scheme on a uniform grid
(finite element, finite difference or finite volume) and applies a
discrete multiresolution decomposition to the numerical data in
order to {\it compress} the computational time and memory space
while preserving the accuracy of the initial scheme. Here the
adaptive sets $\Lambda_n$ are therefore limited within the
resolution level of the uniform grid where the classical scheme
operates. This approach seems more appropriate for hyperbolic
initial value problems \cite{CKMP,DGM}, in which a straightforward
wavelet discretization might fail to converge properly. It should
again be compared to its adaptive mesh refinement counterpart such
as in \cite{BO,BC}.

Let us conclude by saying that despite its theoretical success, in
the sense of achieving for certain classes of problems the optimal
convergence rate with respect to the number of degrees of freedom,
the wavelet-based approach to adaptive numerical simulation
suffers from three major curses. \nl \nl {\bf The curse of
geometry :} while the construction of wavelet bases on a
rectangular domains is fairly simple---one can use tensor product
techniques and inherit the simplicity of the univariate
construction---it is by far less trivial for domains with
complicated geometries. Several approaches have been proposed to
deal with this situation, in particular domain decomposition into
rectangular patches or hierarchical finite element spaces, see
\cite{Co, Da}, and concrete implementations are nowaday available,
but they result in an unavoidable loss of structural simplicity in
comparison to the basic Haar system of \S 3. \nl \nl {\bf The
curse of data structure :} encoding and manipulating the adaptive
wavelet approximations $U^n$ to the solution means that we both
store the coefficients and the indices of the adaptive set
$\Lambda_n$ which should be dynamically updated. The same goes for
the indices of the matrix $A$ which are used in the matrix-vector
algorithm \iref{approx2} at each step of the algorithm. This
dynamical adaptation, which requires appropriate data structure,
results in major overheads in the computational cost which are
observed in practice : the numerical results in \cite{BBCCDDU}
reveal that while the wavelet adaptive algorithm indeed exhibits
the optimal rate of convergence and slightly outperforms adaptive
finite element algorithms from this perspective, the latter
remains significantly more efficient from the point of view of
computational time. \nl \nl {\bf The curse of anisotropy :}
adaptive wavelet approximation has roughly speaking the same
properties as isotropic refinement. However, many instances of
singularities such as boundary layers and shock waves, have
anisotropic features which suggests that the refinement should be
more pronounced in one particular direction. From a theoretical
point of view, the following example illustrate the weakness of
wavelet bases in this situation : if $f=\Chi_\Omega$ with
$\Omega\subset \RR^d$ a smooth domain, then the rate of best
$N$-term approximation in $X=L^2$ is limited to $r=1/(2d-2)$ and
therefore deteriorates as the dimension grows. Wavelet bases
should therefore be reconsidered if one wants to obtain better
rates which take some advantage of the geometric smoothness of the
curves of dicsontinuities. On the adaptive finite element side,
anisotropic refinement has been considered and practically
implemented, yet without a clean theory available for the design
of an optimal mesh. \nl

The significance of wavelets in
numerical analysis remains therefore tied to these curses
and future breakthrough are to be expected once
simple and appropriate solutions are proposed
in order to deal with them.

\label{lastpage}


\begin {thebibliography} {99}

\bibitem{BR} Babushka, I. and W. Reinhbolt (1978) {\it A-posteriori analysis for adaptive
finite element computations}, SIAM J. Numer. Anal. 15, 736--754.

\bibitem{BM}
Babu\v{s}ka, I. and A.  Miller (1987), {\it A feedback finite element
method with a-posteriori error estimation: Part I. The finite element method and some
basic properties of the a-posteriori error estimator}, Comput.\ Methods
Appl.\ Mech.\ Engrg.\ {\bf 61}, 1--40.

 \bibitem{BW} Bank, R.E. and A. Weiser (1985),  {\it Some a posteriori error estimates
for elliptic
partial differential equations}, Math.\ Comp., {\bf 44}, 283--301.

\bibitem{BBCCDDU}
Barinka, A., T. Barsch, P. Charton, A. Cohen, S. Dahlke, W.
Dahmen, K. Urban (1999), {\it Adaptive wavelet schemes for
elliptic problems---Implementation and numerical experiments, IGPM
Report \# 173 RWTH Aachen}, SIAM J. Sci. Comp. {\bf 23}, 910--939.


\bibitem{BC} Berger, M. and P. Collela (1989) {\it Local adaptive mesh
refinement for shock hydrodynamics}, J. Comp. Phys. 82, 64--84.

\bibitem{BO} Berger, M. and J. Oliger (1984) {\it Adaptive mesh refinement for
hyperbolic partial differential equations}, J. Comp. Phys. 53,
482--512.

\bibitem{Ber} Bertoluzza, S. (1995) {\it A posteriori error estimates for wavelet
Galerkin methods}, Appl. Math. Lett. 8, 1--6.

\bibitem{Be} Bertoluzza, S. (1997)
{\it An adaptive collocation method based on interpolating
wavelets}, in {\it Multiscale Wavelet Methods for PDEs},
W.~Dahmen, A.~J.~Kurdila, P.~Oswald (eds.), Academic Press,
109--135.

\bibitem{BH} Bihari, B. and A. Harten (1997) {\it
Multiresolution schemes for the numerical solution of 2--D
conservation laws}, SIAM J. Sci. Comput. {\bf 18}, 315--354.

\bibitem{BDD} Binev, P., W. Dahmen and R. DeVore (2002), {\it Adaptive
finite element methods with convergence rates}, preprint IGPM-RWTH
Aachenm, to appear in Numerische Mathematik.

\bibitem{Ca} Canuto, C. and I. Cravero (1997), {\it Wavelet-based adaptive methods
for advection-diffusion problems}, Math. Mod. Meth. Appl. Sci.
{\bf 7}, 265--289.

\bibitem{Cia} Ciarlet, P.G. (1991), {\it Basic error estimates for the finite element
method}, Handbook of Numerical Analysis, vol II,
P. Ciarlet et J.-L. Lions eds., Elsevier, Amsterdam.

\bibitem{Co} Cohen, A. (2000), {\it Wavelets in numerical analysis},
Handbook of Numerical Analysis, vol. VII, P.G. Ciarlet and J.L.
Lions, eds., to appear in 2003 as a book``Numerical analysis of
wavelet methods", Elsevier, Amsterdam.

\bibitem{CDD1} Cohen, A., W. Dahmen and R. DeVore (2000), {\it Adaptive
wavelet methods for elliptic operator equations---convergence
rate}, Math. Comp. {\bf 70}, 27--75.

\bibitem{CDD2} Cohen, A., W. Dahmen and R. DeVore (2002), {\it Adaptive wavelet
methods for operator equations---beyond the elliptic case}, Found.
of Comp. Math. {\bf 2}, 203--245.

\bibitem{CDD3} Cohen, A., W. Dahmen and R. DeVore (2002), {\it Adaptive wavelet
methods for nonlinear variational problems}, preprint
IGPM-RWTH Aachen, submitted to SIAM J. Num. Anal.

\bibitem{CKMP} Cohen, A., S.M. Kaber, S. Mueller and M. Postel (2002),
{\it Fully adaptive multiresolution finite volume schemes for
conservation laws}, Math. of Comp. {\bf 72}, 183--225.

\bibitem{CM} Cohen, A. and R. Masson (1999),
{\it Wavelet adaptive methods for elliptic
problems---preconditionning and adaptivity}, SIAM J. Sci. Comp.
21, 1006--1026.

\bibitem{DDU} Dahlke, S., W. Dahmen and K. Urban (2001),
{\it Adaptive wavelet methods  for saddle point problems---optimal
convergence rates}, preprint IGPM-Aachen

\bibitem{DD} Dahlke, S. and R. DeVore (1997),
{\it Besov regularity for elliptic boundary value problems},
Communications in PDEs 22, 1--16.

\bibitem{Da} Dahmen, W. (1997), {\it Wavelet and multiscale methods for operator equations},
Acta Numerica 6, 55--228.

\bibitem{DGM} Dahmen, W., B. Gottschlich-M\"uller and S. M\"uller (2001),
{\it Multiresolution Schemes for Conservation Laws}, Numerische
Mathematik 88, 399--443.

\bibitem{Dau} Daubechies, I. (1992), {\it Ten lectures on wavelets}, SIAM, Philadelphia.

\bibitem{De} DeVore, R. (1997), {\it Nonlinear Approximation},
Acta Numerica 51--150.

\bibitem{DeLu} DeVore, R. and B. Lucier (1990), {\it High order regularity for conservation
laws}, Indiana J. Math. 39, 413--430.

\bibitem{Do} D\"orfler, W.  (1996), {\it A convergent adaptive algorithm for Poisson's equation},
SIAM J. Num. Anal. {\bf 33}, 1106--1124

\bibitem{EEHJ} Eriksson, K., D. Estep, P. Hansbo, and C. Johnson (1995),
{\it Introduction to adaptive methods for differential equations}, Acta
Numerica {\bf 4}, Cambridge University Press, 105--158.

\bibitem{H1} Harten, A. (1994),
{\it Adaptive multiresolution schemes for shock computations}, J.
Comp. Phys. 115, 319--338.

\bibitem{Ha} Harten, A. (1995), {\it Multiresolution algorithms for the numerical solution
of hyperbolic conservation laws}, Comm. Pure and Appl. Math. 48,
1305--1342.

\bibitem{MNS} Morin, P., R. Nocetto and K. Siebert (2000), {\it Data oscillation and convergence of
adaptive FEM}, SIAM J. Num. Anal. {\bf 38}, 466--488.

\bibitem{Pet} Petrushev, P. (1988), {\it Direct and converse theorems for spline and
rational approximation and Besov spaces}, in {\it Function spaces
and applications}, M. Cwikel, J. Peetre, Y. Sagher and H. Wallin,
eds., Lecture Notes in Math. 1302, Springer Verlag, Berlin,
363--377.

\bibitem{Tem} Temlyakov, V. (1998), {\it Best $N$-term approximation and greedy algorithms},
Adv. Comp. Math. {\bf 8}, 249--265.

\bibitem{Ve} Verf\"urth, R. (1994), {\it A-posteriori error estimation and adaptive mesh
refinement techniques}, Jour. Comp. Appl. Math. 50, 67--83.

\end{thebibliography}
\end{document}